\documentclass[11pt]{article}
\usepackage{amsmath}
\usepackage{amsfonts}
\usepackage{amssymb}
\usepackage[all]{xy}
\usepackage{stmaryrd}

\newtheorem{definition}{Definition}
\newtheorem{proposition}{Proposition}

\newenvironment{proof}{{\bf Proof:}}{$\text{ }\blacksquare$}

\begin{document}
\title{The club of simplicial sets}
\author{Dennis Borisov\\ dennis.borisov@gmail.com\\ Max-Planck Institute for Mathematics, Bonn, Germany}
\date{\today}
\maketitle

\begin{abstract}
A club structure is defined on the category of simplicial sets. This club generalizes the operad of associative rings by adding ``amalgamated'' products.
\end{abstract}

\section{Introduction}

There is a straightforward way to define operads in the monoidal category $(Cat,\times)$: just apply the standard definition. However, since $Cat$ is not just a category, but a $2$-category, such a definition is of very limitied value. For one thing, the action of symmetric groups through functors would not be the correct action in most applications. Instead, one would need symmetric groups to act by morphisms.

In this paper we make use of another way, the $2$-categorical structure of $Cat$ makes itself felt. In the definition of operads one parameterizes the procedure of taking several points in a set and composing them into one point. Of course sets can be substituted with objects in any other symmetric monoidal category, but the principle remains the same: we compose strings of elements.

\smallskip

Let $M\in Cat$ be a category. A string of objects $\{A_1,\ldots,A_n\}\subseteq M$ is the same as a diagram ${\bf n}\rightarrow M$, where ${\bf n}$ is a discrete category on $n$-objects. Here of course we can take any diagram $D\rightarrow M$, where $D$ is not necessarily discrete, and try to ``compose'' it. This kind of compositions cannot be described in terms of operads. We need the notion of a club instead.

\smallskip

In their full generality clubs were developed by G.M.Kelly in \cite{KG74}, and they are the way to encode associativity of compositions when we compose arbitrary diagrams, and not just the ones parameterized by discrete categories.

Recall that the main axiom of operads is associativity of compositions of operations. Clubs provide a formalization for the same axiom, but in the more general case, where operations can have arbitrary diagrams as inputs, and not just strings. Of course, operads are clubs of a particular kind.

\smallskip

In this paper we define a club structure on the category of simplicial sets $SSet$. This club generalizes the operad of associative rings by adding compositions of elements {\it relative} to other elements. 

Discrete simplicial sets give just an associative product, non-discrete simplicial sets give ``amalgamated'' associative products. An example of such amalgamated products are monoidal globular categories in \cite{Ba98}.

\smallskip

Here is the structure of the paper: in section \ref{SemiDirect} we recall the definition of clubs. The general definition, given in \cite{KG74}, takes place in an arbitrary $2$-category. We do not need this generality, and we consider only clubs in $Cat$. As an example we show that set-theoretic operads are clubs in $Cat$ of a particular kind.

We also use a different notation from \cite{KG74}. A club in $Cat$ is a monoid in the category of diagrams in $Cat$. This category has a very important monoidal product, which is a straightforward generalization to categories of the semidirect product of groups. Therefore we use the symbol $\ltimes$ to designate this monoidal structure.

For any diagram $\underline{D}$ in $Cat$, the functor $\underline{D}\ltimes-:Cat\rightarrow Cat$ is an instance of what is called {\it a familial $2$-functor} (\cite{WM07}, \cite{SR00}).

\smallskip

In section \ref{SSet}, starting from the category $SSet$ of simplicial sets, we define the structure of a club on $\underline{SSet}$, where $\underline{SSet}$ is a diagram in $Cat$, parameterized by $SSet$, with every simplicial set mapped to its category of simplices.

In fact, we define two clubs: one on the entire category of simplicial sets, and another on the subcategory, consisting of injective morphisms. The latter is important in applications, when we have a category $M$, with an $\underline{SSet}$-algebra structure on it, and we want to have an $\underline{SSet}$-algebra structure on the subcategory of mono-morphisms in $M$.

\smallskip

\underline{A note on notation:} When working with sets we use the approach of universes (\cite{SGA4}), in particular we speak of {\it a} small category $Set$ of sets, meaning sets in a given universe. Consequently we have {\it a} small category $SSet$ of simplicial sets.


\section{Semi-direct product and clubs in $Cat$}\label{SemiDirect}

\begin{definition}\label{CategoryCategories} Let $\underline{Cat}$ be the following category:\begin{itemize}
\item An object $\underline{D}\in\underline{Cat}$ is a pair $\{D,R\}$, where $D$ is a small category, and $R:D\rightarrow Cat$ is a functor.
\item A morphism $\underline{F}:\underline{D}\rightarrow\underline{D}'$ in $\underline{Cat}$ is a pair $\{F,\rho\}$ where $F:D\rightarrow D'$ is a functor, and $\rho$ is a natural transformation, making the following diagram commutative:
	$$\xymatrix{D\ar[rrrr]^{F}\ar[rrdd]_R &&&& D'\ar[lldd]^{R'}\\ &\ar@{<-}[rr]^\rho &&& \\ && Cat &&}$$
\end{itemize}\end{definition}
Now we are going to define a monoidal structure $\{\ltimes,\underline{\bf 1}\}$ on $\underline{Cat}$. We start with products of objects. 

\smallskip

Let $\underline{D},\underline{D}'\in\underline{Cat}$ be objects, let $d\in D$ be an object, and let $R(d)\overset{\psi}\rightarrow D'$ be a functor. We define a category $R(d)\ltimes_\psi\underline{D}'$ as follows:\begin{itemize}
\item Objects of $R(d)\ltimes_\psi\underline{D}'$ are pairs $\{a,b\}$, where $a\in R(d)$ and $b\in R'\psi(a)$.
\item Morphisms of $R(d)\ltimes_\psi\underline{D}'$ are pairs $\{\alpha,\beta\}$, where $\alpha:a_1\rightarrow a_2$ is a morphism in $R(d)$, and $\beta:R'\psi(\alpha)(b_1)\rightarrow b_2$ is a morphism in $R'\psi(a_2)$.
\item Composition of $\xymatrix{\{a_1,b_1\}\ar[r]^{\{\alpha_1,\beta_1\}} & \{a_2,b_2\}\ar[r]^{\{\alpha_2,\beta_2\}} & \{a_3,b_3\}}$ is
	$$\{\alpha_2\alpha_1,\beta_2R'\psi(\alpha_2)(\beta_1)\}.$$\end{itemize}
It is easy to check that $R(d)\ltimes_\psi\underline{D}'$ is indeed a category, and similarity between this construction and semi-direct product of groups is obvious. 

Different from the case of groups, we can put together all $R(d)\ltimes_\psi\underline{D}'$'s for all $d\in D$ and all $\psi:R(d)\rightarrow D'$ to get a diagram $\underline{D}\ltimes\underline{D}'\in\underline{Cat}$.

\smallskip

First we describe the parameterizing category. Let $D\ltimes D'$ be the small category, defined as follows:\begin{itemize}
\item Objects of $D\ltimes D'$ are pairs $\{d,\psi_d\}$, where $d\in D$ is an object, and $\psi_d:R(d)\rightarrow D'$ is a functor.
\item Morphisms in $D\ltimes D'$ are pairs $\{f,\phi\}$, where $f:d_1\rightarrow d_2$ is a morphism in $D$, and $\phi$ is a natural transformation, making the following diagram commutative:
	$$\xymatrix{R(d_1)\ar[rrrr]^{R(f)}\ar[rrdd]_{\psi_{d_1}} &&&& R(d_2)\ar[lldd]^{\psi_{d_2}}\\ &\ar[rr]^\phi &&&\\ && D' &&}$$
\end{itemize}
Now we define $R\ltimes R':D\ltimes D'\rightarrow Cat$. As we said above, we would like to collect all $R(d)\ltimes_{\psi_d}\underline{D}'$'s into one diagram, so on objects $R\ltimes R'$ is clear:  
	$$\xymatrix{R\ltimes R':\{d,\psi_d\}\ar@{|->}[r] & R(d)\ltimes_{\psi_d}\underline{D}'.}$$ 
It is straightforward then to define the action of $R\ltimes R'$ on morphisms of $D\ltimes D'$ by composing functors and natural transformations in an obvious way. Here is the explicit description: for a morphism $\{d_1,\psi_{d_1}\}\overset{\{f,\phi\}}\rightarrow\{d_2,\psi_{d_2}\}$ we define the functor 
	$$\xymatrix{R\ltimes R'(\{f,\phi\}):R(d_1)\ltimes_{\psi_{d_1}}\underline{D}'\ar[r] & R(d_2)\ltimes_{\psi_{d_2}}\underline{D}'}$$ 
as follows: let $\{\alpha,\beta\}:\{a_1,b_1\}\rightarrow\{a_2,b_2\}$ be a morphism in $R(d_1)\ltimes_{\psi_{d_1}}\underline{D}'$, then we define $R\ltimes R'(\{f,\phi\})(\{\alpha,\beta\})$ to be
	$$\xymatrix{\{R(f)(a_1),R'(\phi_{a_1})(b_1)\}\ar[rrr]^{\{R(f)(\alpha),R'(\phi_{a_2})(\beta)\}} &&& \{R(f)(a_2),R'(\phi_{a_2})(b_2)\}}$$
Functoriality of this construction is obvious, since all we do here is composing functors and natural transformations.

\begin{definition}\label{SemiDirect} Let $\underline{D},\underline{D}'\in\underline{Cat}$ be objects. We define their {\bf semi-direct product} $\underline{D}\ltimes\underline{D}'$ to be $\{D\ltimes D',R\ltimes R'\}\in\underline{Cat}$ as above.\end{definition} 
Now we describe the unit objects for $\ltimes$. Let ${\bf 1}$ be a discrete category on one object, and let $\underline{\bf 1}\in\underline{Cat}$ consist of ${\bf 1}$, mapped to itself in $Cat$. It is easy to see that for any $\underline{D}\in\underline{Cat}$ we have canonically
	$$\underline{D}\ltimes\underline{\bf 1}\cong\underline{D}\cong\underline{\bf 1}\ltimes\underline{D}.$$
	
We do not prove the following proposition, since it is a consequence of the general result, proved in \cite{KG74}.

\begin{proposition} Let $\underline{Cat}$ be the category of small diagrams in $Cat$ (Definition \ref{CategoryCategories}). The semi-direct product $\ltimes$ (Definition \ref{SemiDirect}) together with $\underline{\bf 1}$ define a monoidal structure on $\underline{Cat}$.\end{proposition}

We would like to note that $\ltimes$ is {\it not} symmetric. This is easy to see from the following simple example: let $D={\bf 1}$, $D'={\bf 2}$ (discrete categories on one and two objects respectively); let $R:D\rightarrow Cat$ be defined by by mapping the only object to ${\bf 2}\in Cat$, and let $R':D'\rightarrow Cat$ be defined by mapping every object to ${\bf 1}\in Cat$. Then
	$$\underline{D}\ltimes\underline{D}'\ncong\underline{D}'\ltimes\underline{D}.$$
	
\begin{definition} {\bf A club} in $Cat$ is a monoid in $(\underline{Cat},\ltimes,\underline{\bf 1})$.\end{definition}

As with most monoids, we will be interested in modules over a club in $Cat$. Given a club $\underline{C}$, it is straightforward to define a $\underline{C}$-module in $\underline{Cat}$, but in practice we would like clubs to act on categories, i.e. objects of $Cat$, rather than $\underline{Cat}$. For that we need a bit of notation.

\smallskip

Let $D$ be a small category. There are is a natural way to associate an object in $\underline{Cat}$ to $D$. Let $\mathbb D:=D\mapsto{\bf 1}$ be the diagram in $Cat$, having $D$ as the parameterizing category, s.t. every object in $D$ is mapped to ${\bf 1}\in Cat$.

Notice that the assignment $D\mapsto\mathbb D$ is a functor from $Cat$ to $\underline{Cat}$, and it is left adjoint to the forgetful functor $\underline{Cat}\rightarrow Cat$, that maps every diagram to its parameterizing category.

\begin{definition}\label{Action} Let $\underline{C}$ be a club in $Cat$, and let $M\in Cat$ be a category. We define $\underline C(M)$ to be the parameterizing category of $\underline{C}\ltimes\mathbb M$.

{\bf A $\underline C$-algebra} is a category $M$, together with a functor $\underline{C}(M)\rightarrow M$, satisfying the usual associativity conditions.\end{definition}

\smallskip

Now we are ready to consider examples.\begin{itemize}
\item[1.] Let $P$ be a set-theoretic non-symmetric operad, i.e. we have a sequence of sets $\{P_n\}_{n\geq 0}$, a chosen element $e\in P_1$, and a sequence of compositions
	$$\gamma_{m_1,\ldots,m_n}: P_n\times P_{m_1}\times\ldots\times P_{m_n}\rightarrow P_{m_1+\ldots+m_n},$$
satisfying the usual conditions of associativity and unitality.

Now we construct a diagram in $Cat$, starting with $P$, and for every $n>0$ a choice of a discrete category ${\bf n}$ having an ordered set of $n$ objects. The parameterizing category $P$ is the discrete category having $\underset{n\geq 0}\coprod P_n$ as the set of objects. Every $p\in P_n\subseteq P$ is mapped to ${\bf n}$. We will denote the resulting diagram by $\underline{P}$.

This construction works for any $\mathbb N$-collection in $Set$, in particular for $P\circ P$, where
	$$(P\circ P)_k=\underset{m_1+\ldots+m_n=k}\coprod P_n\times P_{m_1}\times\ldots\times P_{m_n}.$$
Proof of the following proposition is straightforward.	

\begin{proposition}\label{NonSymmetric}\begin{itemize}
\item[1.] For any collection $P$ in $Set$ we have
	\begin{equation}\label{OperadCorrespondence}\underline{P\circ P}\cong\underline{P}\ltimes\underline{P}.\end{equation}
\item[2.] The correspondence (\ref{OperadCorrespondence}) defines a bijection between the set of operadic compositions $\{\gamma_{m_1,\ldots,m_n}\}$ on $P$, and the set of $\ltimes$-monoidal structures 
	$$\{F,\rho\}:\underline{P}\ltimes\underline{P}\rightarrow\underline{P},\qquad\{I,\iota\}:\underline{\bf 1}\rightarrow\underline{P},$$ 
s.t. $\rho,\iota$ are natural equivalences, preserving the order on ${\bf n}$'s.\end{itemize}\end{proposition}

\smallskip

Let $S$ be a set, and suppose $P$ acts on it. Then it is easy to see how to translate such an action into the structure of a $\underline{P}$-algebra on $S$, considered as a discrete category.

\item[2.] Now let $P$ be an operad in $Set$. Here, in addition to choosing a discrete category ${\bf n}$ on $n$ objects, for each $n\geq 1$ we fix an isomorphism 
	$$\mathbb S_n\cong Aut({\bf n}).$$ 
Then we can define a diagram $\underline{P}\in\underline{Cat}$ as follows: the parameterizing category $P$ has $\underset{n\geq 0}\coprod P_n$ as the set of objects, and $\forall p,q\in P_n\subseteq P$, we put $Hom(p,q)$ to be the set of all $\sigma_n\in\mathbb S_n$, s.t. $\sigma_n(p)=q$; the functor $R:P\rightarrow Cat$ maps every $p\in P_n\subseteq P$ to ${\bf n}$, and every morphism $p\overset\sigma\rightarrow q$ to the corresponding endofunctor on ${\bf n}$.

It is clear that $\underline{P}$ is indeed an object in $\underline{Cat}$, and we can apply the same technique to every $\Sigma$-collection in $Set$. Different from the non-symmetric case, we have that in general $\underline{P}\ltimes\underline{P}\ncong\underline{P\circ P}$. However, we have a natural inclusion $\underline{P}\ltimes\underline{P}\rightarrow\underline{P\circ P}$, and hence we can conclude the following.

\begin{proposition} For any $\Sigma$-collection $P$ in $Set$, there is a bijection between operadic structures on $P$, and $\ltimes$-monoidal structures 
	$$\{F,\rho\}:\underline{P}\ltimes\underline{P}\rightarrow\underline{P},\qquad\{I,\iota\}:\underline{\bf 1}\rightarrow\underline{P},$$ 
s.t. $\rho,\iota$ are natural equivalences, that preserve order on ${\bf n}$'s.\end{proposition}
\begin{proof} The only difference here from the non-symmetric case is the action of symmetric groups. Since $P_n\times P_{m_1}\times\ldots\times P_{m_n}$ carries the action of only $\mathbb S_n\times\mathbb S_{m_1}\times\ldots\times\mathbb S_{m_n}$, and hence in general it is not an $\mathbb S_{m_1+\ldots+m_n}$-set, we have that $\underline{P}\ltimes\underline{P}\ncong\underline{P\circ P}$, and we cannot proceed as in Proposition \ref{NonSymmetric}.

In defining operads one extends $P_n\times P_{m_1}\times\ldots\times P_{m_n}$ by tensoring it with $\mathbb S_{m_1+\ldots+m_n}$ over $\mathbb S_n\times\mathbb S_{m_1}\times\ldots\times\mathbb S_{m_n}$. However, while this makes definition of an operad cleaner, it is not really needed, and it is enough to postulate equivariance only with respect to $\mathbb S_n\times\mathbb S_{m_1}\times\ldots\times\mathbb S_{m_n}$.\end{proof}

\smallskip

Also here it is easy to see how to translate the notion of a $P$-algebra in $Set$ into a $\underline{P}$-algebra in $Cat$.
\end{itemize}


\section{The club of simplicial sets}\label{SSet}

In the previous section we have considered two examples of $\ltimes$-monoids of a special kind. In general a $\ltimes$-monoid is given by an object $\underline{D}$ in $\underline{Cat}$, together with morphisms 
	$$\{F,\rho\}:\underline{D}\ltimes\underline{D}\rightarrow\underline{D},\qquad\{I,\iota\}:\underline{\bf 1}\rightarrow\underline{D}$$ 
in $\underline{Cat}$, satisfying the usual associativity and unit axioms. 

In the case of set-theoretic operads we have required that $\rho$ and $\iota$ are not just natural transformations, but natural {\it equivalences}. This requirement was a consequence of the way we represented operads: all operadic compositions were encoded in the parameterizing category $P$, i.e. operations are represented as objects in $P$.  The functor $P\rightarrow Cat$ was there only to keep track of the arity of these operations.

Now we consider a case where $\rho$ is not required to be invertible. This case is the main example for a ``diagrammatic'' operadic action on categories: here we do not compose strings of objects, but diagrams of objects, and hence categories in the image of $P\rightarrow Cat$ stop being just a bookkeeping device, but carry information of their own. 

The diagrams in question here are given by simplicial sets. We start with defining a procedure that produces a category out of a simplicial set.

\smallskip

Let $SSet$ be a small category of simplicial sets. For any $\mathcal S\in SSet$, $\mathcal S=\{\mathcal S_n\}_{n\geq 0}$, we define a category $S$ as follows:\begin{itemize}
\item The set of objects in $S$ is $\underset{n\geq 0}\coprod\mathcal S_n$.
\item Given two objects $s_m\in\mathcal S_m\subseteq S$, $s_n\in\mathcal S_n\subseteq S$, $Hom(s_m,s_n)$ is the set of all simplicial operators $\mathcal S_m\rightarrow\mathcal S_n$, that map $s_m$ to $s_n$.\end{itemize}
It is clear that for any $\mathcal S\in SSet$, $S$ is a small category. It is also clear that for any morphism $f:\mathcal S\rightarrow\mathcal S'$ in $SSet$ there is a functor $R(f):S\rightarrow S'$, and that this assignment $f\mapsto R(f)$ is functorial. Therefore we have an object $\underline{SSet}:=\{SSet,R\}\in\underline{Cat}$.

\begin{proposition} There is a structure of $\ltimes$-monoid on $\underline{SSet}$.\end{proposition}
\begin{proof} We need to define
	\begin{equation}\label{Diagonal}\{\Delta,\delta\}:\underline{SSet}\ltimes\underline{SSet}\rightarrow\underline{SSet},\end{equation}
where $\Delta:SSet\ltimes SSet\rightarrow SSet$ is a functor, and $\delta$ is a natural transformation
	$$\xymatrix{SSet\ltimes SSet\ar[rrrr]^{\Delta}\ar[rrdd]_{R\ltimes R} &&&& SSet\ar[lldd]^R\\ &\ar@{<-}[rr]^\delta &&& \\ && Cat &&}$$
We start with defining $\Delta$ on objects. Let $\{\mathcal S,\psi\}$ be an object in $SSet\ltimes SSet$, i.e. $\mathcal S$ is a simplicial set, and $\psi:S\rightarrow SSet$ is a functor. The set of objects in the category $S\ltimes_\psi\underline{SSet}$ is graded by $\mathbb Z_{\geq 0}\times\mathbb Z_{\geq 0}$, indeed, objects in $S$ are $\mathbb Z_{\geq 0}$-graded by dimension of simplices, and similarly for categories in the image of $R\circ\psi$.

It is easy to see that the category $S\ltimes_\psi\underline{SSet}$ can be obtained from a bisimplicial set $\mathbb T$ by considering simplices as objects and bisimplicial operators as morphisms. Let $\mathcal T\in SSet$ be the diagonal in $\mathbb T$. We set
	$$\Delta(\{\mathcal S,\psi\}):=\mathcal T,\qquad\delta:T\rightarrow S\ltimes_\psi\underline{SSet},$$
with $\delta$ being given by the diagonal. Thus we have defined (\ref{Diagonal}) on objects.

\smallskip

Let $\{f,\phi\}:\{\mathcal S,\psi\}\rightarrow\{\mathcal S',\psi'\}$ be a morphism in $SSet\ltimes SSet$, where $f:\mathcal S\rightarrow\mathcal S'$ is a map of simplicial sets, and $\phi$ is a natural transformation 
	$$\xymatrix{S\ar[rrrr]^{f}\ar[rrdd]_{\psi} &&&& S'\ar[lldd]^{\psi'}\\ &\ar@{->}[rr]^\phi &&& \\ && SSet &&}$$
The pair $\{f,\phi\}$ induces a functor $S\ltimes_\psi\underline{SSet}\rightarrow S'\ltimes_{\psi'}\underline{SSet}$ that preserves the $\mathbb Z_{\geq 0}\times\mathbb Z_{\geq 0}$-grading, and defines a map of bisimplicial sets $\mathbb T\rightarrow\mathbb T'$. Consequently we get a functor $T\rightarrow T'$ and a corresponding map of simplicial sets $\mathcal T\rightarrow\mathcal T'$. This completes definition of (\ref{Diagonal}).

\smallskip

Now we turn to associativity of $\{\Delta,\delta\}$. This is rather straightforward, essentially it amounts to associativity of taking diagonals in bisimplicial sets.

It remains to define the unit. There is an obvious $\underline{\bf 1}\rightarrow\underline{SSet}$, with ${\bf 1}$ going to the $1$-point simplicial set.\end{proof}

\smallskip

Having the $\ltimes$-monoid structure on $\underline{SSet}$, we can talk about $\underline{SSet}$-algebras in $Cat$. For example, if a category $M$ is closed with respect to taking colimits, we have a canonical structure of a $\underline{SSet}$-algebra on $M$: given a simplicial diagram in $M$, take its colimit.

Sometimes, having a $\underline{SSet}$-algebra $M$ we might like to work only with the subcategory of $M$, consisting of mono-morphisms. It might happen that the action of the entire $\underline{SSet}$ does not preserve the chosen subcategory. In these cases the following definition is useful.

\smallskip

Let $M$ be a category, and let $\mathcal I$ be a set of generators of $M$. For example, if $M=SSet$, $\mathcal I$ is the set of standard simplices $\{\Delta[n]\}_{n\geq 0}$. Let $\mathcal D_\mathcal S:\mathcal S\rightarrow\mathcal M$ be an object in $\underline{SSet}(M)$ (Definition \ref{Action}). 

\begin{definition} Let $I\in\mathcal I$. {\bf An $I$-point} in $\mathcal D_\mathcal S$, is a morphism $I[n]\rightarrow\mathcal D_\mathcal S$ in $\underline{SSet}\ltimes M$, where $I[n]$ is $I$, considered as a constant diagram over $\Delta[n]$.\end{definition} 
It is clear that for each such $I\in\mathcal I$ we obtain a simplicial set $I(\mathcal D_\mathcal S)$ of $I$-points.  If $(F,\phi):\mathcal D_\mathcal S\rightarrow\mathcal D_\mathcal T$ is a morphism in $\mathfrak S(\mathcal M)$, it induces a morphism of simplicial sets $(F,\phi)_I:I(\mathcal D_\mathcal S)\rightarrow I(\mathcal D_\mathcal T)$. 

\begin{definition}\label{Fibrations} We will say that $(F,\phi)$ is {\bf a fibration}, if $F$ is injective, and for each $I\in\mathcal I$ the morphism of simplicial sets $(F,\phi)_I$ is a fibration.\end{definition}
It is straightforward to check that $\{I[n]\}_{I\in\mathcal I,n\geq 0}$ is a set of generators for $\underline{SSet}(M)$, and hence we can iterate this definition to get the notion of a fibration in $\underline{SSet}^k(M)$ for any $k\geq 1$. 

Let $sset\subset SSet$ be the subcategory consisting of injective morphisms, and let $\underline{sset}$ be the corresponding object in $\underline{Cat}$. Let $\underline{sset}\circ\underline{sset}\subset\underline{sset}\ltimes\underline{sset}$ to be the subcategory of fibrations. Similarly, we define $\underline{sset}^{\circ^k}$ for any $k\geq 1$.

\begin{proposition} The sequence $\{\underline{sset}^{\circ^k}\}_{k\geq 1}$ is stable with respect to the $\ltimes$-monoid structure on $\underline{SSet}$.\end{proposition}

Using this proposition we can regard $\underline{sset}$ itself as a monoid, and hence consider $\underline{sset}$-algebras. Usually it happens that if $M$ is an $\underline{SSet}$-algebra, then the category of mono-morphisms in $M$ is an $\underline{sset}$-algebra.

\end{document}